\documentclass[12pt,a4paper,reqno]{amsart}

\usepackage[a4paper,margin=24mm]{geometry}
\usepackage{microtype}
\usepackage{amsthm,amsmath,amsfonts,amssymb,mathtools,dsfont}
\usepackage[colorlinks,citecolor=blue,urlcolor=blue,linkcolor=blue]{hyperref}

\theoremstyle{plain}
\newtheorem{theorem}{Theorem}
\newtheorem{lemma}{Lemma}
\newtheorem{proposition}{Proposition}

  \newtheorem{claim}{Claim}

\theoremstyle{definition}

\newcommand{\Z}{\mathbb Z}
\newcommand{\N}{\mathbb N}
\newcommand{\R}{\mathbb R}
\newcommand{\T}{\mathbb T}
\newcommand{\Pp}{\mathbb P}
\newcommand{\E}{\mathbb E}
\newcommand{\ind}[1]{\mathds{1}_{\{#1\}}}
\DeclareMathOperator{\Var}{Var}
\DeclareMathOperator{\Cov}{Cov}

\title[No roughening transition]{No roughening transition for the Ising model on $\mathbb{Z}^d$ with $d \geq 4$}

\author{Jianping Jiang}
\address{Yau Mathematical Sciences Center, Tsinghua University, Beijing 100084, China.}
\email{jianpingjiang@tsinghua.edu.cn}
\author{Sike Lang}
\address{Qiuzhen College, Tsinghua University, Beijing 100084, China.}
\email{langsk24@mails.tsinghua.edu.cn}

\begin{document}
	
\begin{abstract}
Consider the Ising measure $\mu_{\Lambda_n,\beta}^{\mathrm{Dob}}$ on $\Lambda_n:=([-n,n]^{d-1}\times[-n,n-1]) \cap \mathbb{Z}^d$ at inverse temperature $\beta$ with all plus boundary condition at heights $z \geq 0$ and all minus boundary condition at heights $z<0$. Let $\beta_c(d)$ be the critical inverse temperature. For every $d \geq 4$ and $\beta>\beta_c(d)$, we prove that $\mu_{\Lambda_n,\beta}^{\mathrm{Dob}}(\sigma_0)$ has a positive lower bound independent of $n$.  As a consequence, any Dobrushin state (i.e., a subsequential limit of $\mu_{\Lambda_n,\beta}^{\mathrm{Dob}}$) is not translation invariant for those $(d,\beta)$, and thus there is no roughening transition when $d \geq 4$.
\end{abstract}

\maketitle

\section{Introduction and main result}
The Ising model on a finite set $\Lambda\subset\Z^d$, at inverse temperature $\beta\ge0$, with boundary condition $\eta\in\{-1,+1\}^{\Z^d}$ and external field $h\in \R^{\Lambda}$, is the probability measure on $\sigma\in\{-1,+1\}^{\Lambda}$ defined by
\begin{equation}\label{eq:Isingdef}
	\mu^{\eta}_{\Lambda,\beta,h}(\sigma)
	:=\frac{1}{Z^{\eta}_{\Lambda,\beta,h}}
	\exp\Bigg[
	\beta\sum_{\substack{u, v \in \Lambda \\ u \sim v }}\sigma_u\sigma_v
	+\beta\sum_{\substack{u\in\Lambda,\ v\notin\Lambda\\u\sim v}}
	\sigma_u\eta_v
	+\sum_{u\in\Lambda} h_u\sigma_u
	\Bigg],
\end{equation}
where $u\sim v$ means that $u$ and $v$ are nearest neighbors (each edge is counted once), and $Z^{\eta}_{\Lambda,\beta,h}$ is the partition function.  Write $\langle \cdot \rangle ^{\eta}_{\Lambda,\beta,h}$ for the expectation with respect to $\mu^{\eta}_{\Lambda,\beta,h}$. Let $\beta_c(d)$ be the critical inverse temperature of the zero-field Ising model. 

The Dobrushin boundary condition is defined
\begin{equation}\label{eq:Dobdef}
	\eta_u:=\begin{cases}
		+1,& u_d \geq 0,\\
		-1, & u_d <0,
	\end{cases}
\end{equation}
where $u_d$ is the $d$-th coordinate of $u$. For each $n\in \N$, define
\[\Lambda_n:=([-n,n]^{d-1}\times[-n,n-1]) \cap \mathbb{Z}^d\]
and write $\mu_{\Lambda_n,\beta}^{\mathrm{Dob}}$ for the measure defined in \eqref{eq:Isingdef} with boundary condition \eqref{eq:Dobdef} and zero external field $h=0$.

\begin{theorem}\label{thm:main}
	Let $d\geq 4$ and $\beta>\beta_c(d)$. There exists $c=c(d,\beta)>0$ such that
	\[\langle \sigma_0 \rangle^{\mathrm{Dob}}_{\Lambda_n,\beta} \geq c,\qquad \forall n \in \mathbb{N}.\]
	In particular, every subsequential local limit of $\mu_{\Lambda_n,\beta}^{\mathrm{Dob}}$ as $n\rightarrow\infty$ is a Gibbs state which is not translation invariant.
\end{theorem}

In the paper, local convergence is the same as weak convergence in the product topology. The states arising from subsequential local limits of $\mu_{\Lambda_n,\beta}^{\mathrm{Dob}}$ are usually called \textit{Dobrushin states}. Dobrushin~\cite{Dob72} first constructed such states and proved that they are not translation invariant when $d \geq 3$ and $\beta$ is sufficiently large. van Beijeren~\cite{vanBei75} gave a simpler proof for all $d\geq 3$ and $\beta >\beta_c(d-1)$ by correlation inequalities. For $d=2$,  Abraham and Reed~\cite{AR76} (see also Gallavotti~\cite{Gal72}) proved that Dobrushin states are translation invariant for all $\beta \geq 0$; later, Aizenman~\cite{Aiz80} and Higuchi~\cite{Hig81} proved that all Gibbs states on $\Z^2$ are translation invariant. When $d=3$, it is expected \cite{WGL73,BFL82} that there is a \textit{roughening inverse temperature} $\beta_R > \beta_c(3)$ such that Dobrushin states  are translation invariant if $\beta < \beta_R$ but not if $\beta > \beta_R$. It is known that the Gibbs state is unique whenever $\beta\in[0,\beta_c(d)]$ for all $d\geq 2$~\cite{ABF87,Yan52,AF86,ADCS15}; Theorem~\ref{thm:main} thus implies that the roughening transition can only occur when $d=3$.

The existence of non-translation-invariant Gibbs states is closely related to the behavior of the interface induced by the Dobrushin boundary condition. We refer to \cite[Section~3.10.7]{FV18} for the rigorous definition of this interface. For $d=2$ and $\beta>\beta_c(2)$, Greenberg and Ioffe~\cite{GI05} (see also Higuchi~\cite{Hig79}) proved that the interface converges weakly under diffusive scaling to a Brownian bridge. For $d\geq 3$, Gheissari and Lubetzky~\cite{GL22} proved a law of large numbers for the maximum height of the interface in the infinite cylinder for all sufficiently large $\beta$; see also~\cite{GL21,GL23a,GL23b} for further developments. For the three-dimensional Ising model on tilted tori with both slope components positive, Laslier and Lubetzky~\cite{LL26} recently proved a Gaussian free field scaling limit for the centered interface height at sufficiently large $\beta$. 

A related geometric question concerns the coexistence of an infinite plus cluster and an infinite minus cluster. The coexistence was proved in~\cite{JL26b} for all $\beta \in [0,\beta_c(d)+\epsilon_d)$ with $d\geq 8$ and some $\epsilon_d>0$, and in~\cite{JL26a} for all $\beta$ less than the critical inverse temperature for triangular-type 3D lattices.

Throughout the paper, we write $D:=d-1$ and identify $\Z^d$ with $\Z^D \times \Z$. We call the first $D$ coordinates \textit{horizontal coordinates}, and the last coordinate the \textit{vertical coordinate}.

The main idea of the proof is to compare a Dobrushin state on $\Z^d$ with another Gibbs state $\nu$ on the same space such that the magnetization difference between the origin $(0,0)$ and $(0,-1)$ of the former dominates that of $\nu$ while $\nu(\sigma_{(0,0)}-\sigma_{(0,-1)})>0$. This is motivated by van Beijeren's proof~\cite{vanBei75}, which essentially compares the Dobrushin state with the plus state on $\Z^{d-1}$. To construct $\nu$, we first define Ising measures $\nu_n$ on the slab $\Z^{d-1}\times [-n,n-1]$ with antiperiodic boundary condition in the vertical direction. Grouping the spins in the vertical direction into a vector, we are able to apply the infrared bound of Fr\"ohlich, Simon and Spencer~\cite{FSS76} to show that the variance of the total magnetization along one vertical line segment grows at most linearly in its size. Meanwhile, $\nu_n$ satisfies one-sided bounds on reflected partial sums (see \eqref{eq:sumnu_nmag} below). With these ingredients, we finally apply Bodineau's classification of translation-invariant Gibbs states~\cite{Bod06} to conclude that the magnetization difference for $\nu$ is positive.

In Section~\ref{sec:slabandcov}, we first apply the infrared bound to derive a covariance estimate for general spin systems on $D$-dimensional tori (here the assumption that $D\geq3$ is essential), and then use this and convexity of the free energy to construct the Ising slab measure $\nu_n$. In Section~\ref{sec:correlation}, we use the positive association of the agreement set of reflected spins to derive correlation inequalities which enable us to compare Dobrushin states and $\nu_n$. Finally, in Section~\ref{sec:proof}, we prove our main result.

\section{A Gibbs state on an antiperiodic slab and a covariance estimate}\label{sec:slabandcov}
Fix $D \geq 3$ and $\beta>0$. In this section, we apply the infrared bound~\cite{FSS76} to columns of spins in an antiperiodic slab. We first derive a covariance bound for general spin measures on $\Z^D$ arising from local limits of periodic measures on tori. Let $\rho$ be a nonzero finite measure supported on a compact subset $\Omega \subset \R^N$ for some $N\in \N$. Let $\T_L:=(\Z / L\Z)^D$ be the torus with length $L \geq 3$. For $h\in \R^N$, consider the spin measure on $\Omega^{\T_L}$,
\begin{equation}\label{eq:P_Lhdef}
	\Pp_{L,h}(dS) = \frac{1}{Z_{L,h}}\exp\Bigg[\beta\sum_{\substack{x,y\in \T_L: x\sim y}} S_x \cdot S_y + \sum_{x\in \T_L} h\cdot S_x\Bigg]\prod_{x\in \T_L} \rho(d S_x),
\end{equation}
where $Z_{L,h}$ is the partition function. Let $\E_{L,h}$ be the expectation with respect to $\Pp_{L,h}$. 

Define the finite-volume \textit{free energy}
\[f_L(h):=\frac{1}{L^D} \ln Z_{L,h}.\]
It is well-known that the infinite-volume \textit{free energy}
\begin{equation}\label{eq:inffreeenergy}
	f(h):=\lim_{L\rightarrow\infty}f_L(h)
\end{equation}
exists and is convex in $h \in \R^N$; see, e.g., Theorem 3.6 of \cite{FV18}, whose proof can be extended to our case. Therefore, $f$ is differentiable at Lebesgue-almost every point. Define
\[G_D(0):=\frac{1}{(2\pi)^D}\int_{[-\pi,\pi]^D}\frac{dk}{2D-2\sum_{j=1}^D \cos(k_j)},\]
which is finite if $D\geq 3$.

\begin{proposition}\label{prop:infrared}
	Suppose the free energy $f$ is differentiable at $h$. Then every subsequential local limit of $\Pp_{L,h}$ as $L\rightarrow\infty$, denoted by $\Pp_h$ (with expectation $\E_h$), is a translation invariant Gibbs state and satisfies
	\begin{equation}\label{eq:varine}
		\E_h(S_0)=\nabla f(h), \Var_{\Pp_h}(a\cdot S_0) \leq G_D(0)\|a\|_2^2/\beta, \qquad \forall a \in \R^N.
	\end{equation}
	Let $h_j$ be a sequence of differentiability points of $f$ with $h_j\rightarrow 0$ as $j\rightarrow\infty$. Then every subsequential local limit of $\Pp_{h_j}$ as $j\rightarrow\infty$ is a Gibbs state at $h=0$ and satisfies the inequality in \eqref{eq:varine}.
\end{proposition}
\begin{proof}
	We want to apply~\cite[Theorem 2.1]{FSS76} with single-site base measure $\exp[h\cdot S] \rho(dS)$. Recall the definitions of discrete partial derivatives and Laplacian from~\cite{FSS76}: for $g: \T_L \rightarrow \R^N$,
	\[(\partial^L_i g)(x):=g(x+e_i)-g(x),\quad (-\Delta^L g)(x):=2Dg(x)-\sum_{i=1}^D[g(x+e_i)+g(x-e_i)],\  \forall x \in \T_L,\]
	where $e_1,\dots,e_D$ are the standard basis vectors for $\Z^D$. For $b=(b_x)_{x\in \T_L}\in (\R^N)^{\T_L}$ with $\sum_{x\in \T_L} b_x=0$, let $p: \T_L \rightarrow \R^N$ be the unique solution of the equation
	\[-\Delta^L p=b \text{ with }\sum_{x\in \T_L}p(x)=0.\]

	Define $g_i(x):=p(x-e_i)-p(x)$ for each $x\in \T_L$. Then $\sum_{i=1}^D \partial_i^L g_i(x)=b(x):=b_x$ and
	\begin{equation*}
		\sum_{x\in \T_L}\sum_{i=1}^D\|g_i(x)\|_2^2=\sum_{x,i} g_i(x)\cdot[p(x-e_i)-p(x)]=\sum_{x\in \T_L} p(x)\cdot b(x)=\langle b, p \rangle=\langle b, (-\Delta^L)^{-1}b \rangle,
	\end{equation*}
	where the inner product is defined in $\ell^2(\T_L;\R^N)$. Then~\cite[Theorem 2.1]{FSS76} implies
	\begin{equation}\label{eq:infraredbd}
		\E_{L,h}\exp\left[\sum_{x\in \T_L} b_x \cdot S_x\right]\leq \exp\left[\frac{1}{2\beta}\langle b,(-\Delta^L)^{-1}b \rangle \right].
	\end{equation}
	
	Let $\T_L^*:=\{2\pi k/L: k\in \{0,1,\dots,L-1\}^D\}$. Define the discrete Green function
	\[G^L(x):=\frac{1}{L^D}\sum_{k\in \T_L^*\setminus\{0\}}\frac{e^{ik\cdot x}}{2D-2\sum_{j=1}^D \cos(k_j)},\qquad \forall x\in \T_L.\]
	Then we have $(-\Delta^L G^L)(x)=\ind{x=0}-1/L^D$ and $\sum_{x\in T_L} G^L(x)=0$. 
	
	For $a \in \R^N$ and $t\in \R$, let $b_x:=ta(\ind{x=0}-L^{-D})$. We have
	\[\langle b, (-\Delta^L)^{-1}b \rangle=\sum_{x\in \T_L} b_x \cdot (taG^L(x))=\sum_{x\in \T_L} ta(\ind{x=0}-L^{-D}) \cdot (taG^L(x))=t^2\|a\|_2^2 G^L(0).\]
	
	Let $\overline{S}_L:=\sum_{x\in \T_L} S_x/L^D$. Then, \eqref{eq:infraredbd} becomes
	\[\E_{L,h}\exp[ta\cdot (S_0-\overline{S}_L)]\leq \exp[(2\beta)^{-1}t^2\|a\|_2^2 G^L(0)].\]
	Translation invariance gives $\E_{L,h}[a\cdot (S_0-\overline{S}_L)]=0$. Taylor expanding both sides of the last displayed inequality to second order in $t$ gives
	\[\E_{L,h}\left[(a \cdot (S_0-\overline{S}_L))^2\right] \leq \beta^{-1}\|a\|_2^2 G^L(0).\]
	
	The Riemann sums
	\[G^L(0)=\frac{1}{(2\pi)^D}\left(\frac{2\pi}{L}\right)^D \sum_{k\in \T_L^*\setminus\{0\}}\frac{1}{2D-2\sum_{j=1}^D \cos(k_j)}\rightarrow G_D(0) \text{ as } L\rightarrow \infty.\]
	Indeed, writing each $k_j$ modulo $2\pi$ in $[-\pi,\pi)$, we can bound the contribution from $0<|k|\leq \delta$ by $C_D \delta^{D-2}$, uniformly in $L$. Therefore,
	\begin{equation}\label{eq:atimesS_0minusmean}
		\limsup_{L\rightarrow\infty}\E_{L,h}\left[(a \cdot (S_0-\overline{S}_L))^2\right] \leq \beta^{-1}\|a\|_2^2 G_D(0).
	\end{equation}
	
	Let $m:=\nabla f(h)$. For $a \in \R^N$, $\delta>0$ and $t>0$, the exponential Markov inequality implies 
	\[\limsup_{L\rightarrow \infty}L^{-D} \ln \Pp_{L,h}(a\cdot \overline{S}_L \geq a\cdot m+\delta) \leq f(h+ta)-f(h)-t(a\cdot m+\delta).\]
	Since $f$ is differentiable at $h$,
	\[f(h+ta)-f(h)-t(a\cdot m+\delta) \leq -t\delta/2, \qquad \text{for all small }t>0.\]
	Choosing a small $t_0>0$ which may depend on $a$ and $\delta$, we have
	\[f_L(h+t_0a)-f_L(h)-t_0(a\cdot m+\delta) \leq -t_0\delta/4, \qquad \text{for all sufficiently large }L.\]
	Therefore,
	\[\Pp_{L,h}\left(a\cdot (\overline{S}_L-m) \geq \delta\right) \leq \exp[-t_0\delta L^D/4], \qquad \text{for all sufficiently large }L.\]
	With $a$ replaced by $-a$, we get
		\[\Pp_{L,h}\left(a\cdot (\overline{S}_L-m) \leq -\delta\right) \leq \exp[-t_0\delta L^D/4], \qquad \text{for all sufficiently large }L.\]
	
	So for each $\epsilon>0$, we have (for a modified $t_0>0$)
	\[\Pp_{L,h}\left(\|\overline{S}_L-m\|_2 \geq \epsilon \right)\leq \sum_{j=1}^N \Pp_{L,h}\left(|(\overline{S}_L-m)_j| \geq \epsilon/\sqrt{N} \right) \leq 2N \exp[-t_0\epsilon L^D/(4\sqrt{N})].\]
	This implies
	\begin{equation}\label{eq:2ndmoment0}
		\lim_{L\rightarrow\infty} \E_{L,h}\|\overline{S}_L-m\|_2^2=0.
	\end{equation}
	
	By translation invariance,
	\begin{equation}\label{eq:S_0mean}
		\E_{L,h}(S_0)=\E_{L,h}(\overline{S}_L),\qquad \Cov_{\Pp_{L,h}}(a\cdot S_0, a\cdot \overline{S}_L)=\Var_{\Pp_{L,h}}(a\cdot \overline{S}_L).
	\end{equation}
	Consequently,
	\begin{equation}\label{eq:varaS_0}
		\Var_{\Pp_{L,h}}(a\cdot S_0)=\E_{L,h}\left[(a \cdot (S_0-\overline{S}_L))^2\right]+\Var_{\Pp_{L,h}}(a \cdot \overline{S}_L).
	\end{equation}
	
	By the Cauchy-Schwarz inequality and \eqref{eq:2ndmoment0},
	\begin{align*}
		\Var_{\Pp_{L,h}}(a \cdot \overline{S}_L)&=\Var_{\Pp_{L,h}}\left(a \cdot (\overline{S}_L-m)\right)\leq \E_{L,h}\left[(a \cdot (\overline{S}_L-m))^2\right] \\
		&\leq \|a\|_2^2 \E_{L,h}\|\overline{S}_L -m\|_2^2 \rightarrow 0 \text{ as }L \rightarrow \infty.
	\end{align*}
Combining this with \eqref{eq:atimesS_0minusmean} and \eqref{eq:varaS_0}, we get
\[\limsup_{L\rightarrow \infty}\Var_{\Pp_{L,h}}(a\cdot S_0) \leq \beta^{-1}\|a\|_2^2 G_D(0).\]

This proves the inequality in \eqref{eq:varine}. The equality in \eqref{eq:varine} follows from
\[\E_{L,h}(S_0)=\nabla f_L(h) \rightarrow \nabla f(h) \text{ as }L\rightarrow \infty\]
and local convergence. Compactness gives a subsequential limit along $h_j\rightarrow 0$; see, e.g., \cite[Section~6.10.2]{FV18}. Continuity of the finite-range specifications implies each such limit satisfies the DLR equations and thus is a Gibbs state at $h=0$. The local limit gives the convergence of the first and second moments of $a\cdot S_0$; thus the inequality in \eqref{eq:varine} survives in the limit along $h_j\rightarrow 0$.
\end{proof}

For each integer $n \geq 2$, let $I_n:=\{-n,\dots,n-1\}$. For each horizontal vertex $x\in \Z^D$, group the spins into a column $S_x:=(\sigma_{x,z})_{z\in I_n} \in \{-1,+1\}^{I_n}$. Define the single-site measure
\begin{equation}\label{eq:rho_ndef}
	\rho_n(S_x)=\exp\left[\beta \left(\sum_{z=-n}^{n-2}\sigma_{x,z}\sigma_{x,z+1} -\sigma_{x,n-1}\sigma_{x,-n}\right)\right].
\end{equation}
Together with the horizontal interactions $\beta\sum_{xy}S_x \cdot S_y$, this defines the Ising model on $\T_L \times I_n$ with antiperiodic boundary condition in the vertical direction. Let $\mu_{L,n,h}$ be the probability measure defined in \eqref{eq:P_Lhdef} with single-site base measure \eqref{eq:rho_ndef}; and let $f^I_n$ be its infinite-volume free energy defined as in \eqref{eq:inffreeenergy} where $I$ stands for Ising. Define the antiperiodic shift $U$ by
\begin{equation*}
	(US_x)_z:=\begin{cases}
		\sigma_{x,z-1},& -n<z\leq n-1,\\
		-\sigma_{x,n-1},&z=-n,
	\end{cases}
	\qquad \forall x \in \T_L.
\end{equation*}

Since $(US_x)\cdot (US_y)=S_x \cdot S_y$ for any $S_x,S_y\in \{-1,+1\}^{I_n}$, $U$ is an orthogonal transformation and its transpose $U^T$ is equal to $U^{-1}$:
\begin{equation*}
	(U^TS_x)_z:=\begin{cases}
		\sigma_{x,z+1},& -n \leq z< n-1,\\
		-\sigma_{x,-n},&z=n-1,
	\end{cases}
	\qquad \forall x \in \T_L.
\end{equation*}
Note that for each $k\in \Z$, $S_x \mapsto U^k S_x$ for each $x\in \T_L$ preserves the measure $\rho_n$ in \eqref{eq:rho_ndef} and $h\cdot (U^k S_x)=((U^k)^Th)\cdot S_x$ for each $S_x \in \{-1,+1\}^{I_n}$ and $h\in \R^{I_n}$. Changing variables gives $Z_{L,n,h}=Z_{L,n,(U^k)^Th}$, and passing to the limit gives
\begin{equation}\label{eq:freeenergyequal}
	f_n^I\left((U^k)^T h\right)=f^I_n(h),\qquad \forall k\in \Z, h\in \R^{I_n}.
\end{equation}

Define
\begin{equation*}
	\tau_z:=\begin{cases}
		1, & z\geq 0, z\in\Z,\\
		-1,& z<0, z\in\Z.
	\end{cases}
\end{equation*}
For each $j \in \N$, choose a differentiability point $h_j$ of $f_n^I$ such that
\begin{equation}\label{eq:h_jdef}
	\frac{\tau_z}{j}< h_j(z) <\frac{\tau_z}{j}+\frac{1}{j^2}, \qquad \forall z \in I_n.
\end{equation} 
Such $h_j$ exist by almost everywhere differentiability of $f_n^I$. Then we have
\begin{equation}\label{eq:h_jproperties}
	jh_j\rightarrow \tau \text{ as } j\rightarrow \infty,\ h_j(z)+h_j(-1-z)>0, \ h_j(z)-h_j(-1-z)>0, \qquad \forall 0\leq z<n.
\end{equation}

For each $h_j$, let $\nu_{n,h_j}$ be a subsequential local limit of $\mu_{L,n,h_j}$ as $L\rightarrow \infty$. Then take a  subsequential local limit of $\nu_{n,h_j}$ as $j\rightarrow\infty$, denoted by $\nu_n$. By the same argument as in the proof of Proposition \ref{prop:infrared}, the state $\nu_n$ is a horizontally translation invariant Gibbs state for the antiperiodic slab at zero external field.

\begin{proposition}\label{prop:covestimate}
	For each integer $n\geq 2$, $\nu_n$ satisfies
	\begin{equation*}
		 \Var_{\nu_n}\left(\sum_{z \in I_n} a_z \sigma_{0,z}\right) \leq G_D(0)\sum_{z\in I_n} a_z^2/\beta, \qquad \forall a \in \R^{I_n},
	\end{equation*}
	and 
	\begin{equation}\label{eq:sumnu_nmag}
		\sum_{z=-M}^{-1} \nu_n(\sigma_{0,z}) \leq 0 \leq \sum_{z=0}^{M-1} \nu_n(\sigma_{0,z}), \qquad 1\leq M \leq n.
	\end{equation}
\end{proposition}
\begin{proof}
	The variance estimate follows from Proposition~\ref{prop:infrared} directly. Let
	\[m_j:=\nu_{n,h_j}(S_0)=\nabla f_n^I(h_j).\]
	Convexity and \eqref{eq:freeenergyequal} imply
	\begin{align*}
		0=f_n^I\left((U^k)^T h_j\right)-f^I_n(h_j)&\geq \left(\nabla f_n^I(h_j)\right)\cdot \left((U^k)^Th_j-h_j\right)=m_j \cdot \left(U^{-k}h_j-h_j\right)\\
		&=h_j \cdot \left(U^k m_j-m_j\right), \qquad \forall k\in \Z.
	\end{align*}
	
	Multiplying by $j$ and passing to the subsequence defining $\nu_n$, we obtain
	\begin{equation}\label{eq:taudotmequal}
		\tau \cdot m \geq \tau \cdot (U^k m),\ m:=\nu_n(S_0), \qquad \forall k\in \Z.
	\end{equation}
	For $1\leq M \leq n$, a direct computation gives
	\[\tau \cdot (U^Mm-m)=2\sum_{z=-M}^{-1} m_z,\  \tau \cdot (U^{-M}m-m)=-2\sum_{z=0}^{M-1} m_z.\]
	Together with \eqref{eq:taudotmequal}, these identities prove \eqref{eq:sumnu_nmag}.
\end{proof}

\section{Correlation inequalities}\label{sec:correlation}
Fix $D\geq 1$ and $\beta \geq 0$. The method in this section is motivated by \cite{vanBei75} and \cite{Leb74}. We first consider the Ising measure $\mu_{L,n,h}$ on $\T_L \times I_n$ with site-dependent external fields such that $\mu_{L,n,h}(\sigma)$ is proportional to
\begin{equation}\label{eq:mu_Lnhdef}
	\exp\Bigg[\beta \sum_{x\in \T_L} \left(\sum_{z=-n}^{n-2}\sigma_{x,z}\sigma_{x,z+1} -\sigma_{x,n-1}\sigma_{x,-n}\right)+\beta \sum_{z=-n}^{n-1} \sum_{\substack{x,y \in \T_L\\ x\sim y}}\sigma_{x,z}\sigma_{y,z}+\sum_{u\in \T_L\times I_n} h_u\sigma_u\Bigg].
\end{equation}
This measure is more general than the Ising measure defined after \eqref{eq:rho_ndef} (though with the same notation) since the external field there is independent of the horizontal coordinates. For $u=(x,z)$ with $z\geq 0$, let $\bar{u}:=(x,-1-z)$ be the reflected image of $u$ with respect to the hyperplane $u_d=-1/2$. Define
\[s_u:=\frac{\sigma_u+\sigma_{\bar{u}}}{2}, \ t_u:=\frac{\sigma_u-\sigma_{\bar{u}}}{2}, \ H_u:=h_u+h_{\bar{u}},\ K_u:=h_u-h_{\bar{u}}.\]

Note that $(s_u,t_u)$ takes a value in $\Omega_0:=\{(\pm 1,0),(0,\pm 1)\}$. Let $V_+:=\T_L \times \{0,\dots,n-1\}$. Under this change of variables, we obtain a probability measure on $\Omega_0^{V_+}$ such that $\mu_{L,n,h}(s,t)$ is proportional to
\[\exp\Bigg[2\beta \sum_{\substack{u,v\in V_+ \\ u\sim v}}(s_us_v+t_ut_v)+\beta \sum_{\substack{u=(x,z)\\ u\in V_+}}(\ind{z=0}-\ind{z=n-1})(s_u^2-t_u^2)+\sum_{u\in V_+}(H_us_u+K_ut_u)\Bigg].\]
We abuse the notation $\mu_{L,n,h}$ for both the Ising measure and the new measure on $\Omega_0^{V_+}$; this should cause no confusion since the arguments specify which measure we use. Let $\langle \cdot \rangle_{L,n,h}$ be the expectation with respect to $\mu_{L,n,h}$.
\begin{lemma}\label{lem:mon}
	Suppose the external field $h\in \R^{\T_L \times I_n}$ satisfies $H_u \geq 0$ and $K_u \geq 0$ for each $u \in V_+$. Then
	\[\langle t_u \rangle_{L,n,h}\geq 0, \  \frac{\partial}{\partial K_v}\langle t_u \rangle_{L,n,h}\geq 0, \  \frac{\partial}{\partial H_v}\langle t_u \rangle_{L,n,h}\leq 0, \qquad \forall u, v\in V_+.\]
\end{lemma}
\begin{proof}
	Define the agreement set
	\[A:=\{u\in V_+: \sigma_u=\sigma_{\bar{u}}\}=\{u\in V_+: t_u=0\}.\]
	Then for any real-valued functions $\phi$ and $\psi$, we have
	\begin{equation}\label{eq:stexpectation}
		\langle \phi(s) \psi(t) \rangle_{L,n,h}=\sum_{A\subset V_+} P(A) \langle \phi(s) \rangle_{A,2\beta,H} \langle \psi(t) \rangle_{A^c,2\beta,K},
	\end{equation}
	where
	\[P(A):=\frac{\prod_{u\in V_+}w_u(A) Z_{A,2\beta, H} Z_{A^c,2\beta, K}}{Z_{L,n,h}}\]
	with
	\begin{equation*}
		w_{x,z}(A):=\begin{cases}
			e^{\beta} \ind{(x,z)\in A}+e^{-\beta} \ind{(x,z) \notin A},& z=0,\\
			e^{-\beta} \ind{(x,z)\in A}+e^{\beta} \ind{(x,z) \notin A},& z=n-1,\\
			1, & 0<z<n-1,
		\end{cases}
		\qquad \forall (x,z)\in V_+.
	\end{equation*}
	
	Here $\langle \cdot \rangle_{A,2\beta, H}$ is the Ising measure on the subgraph induced by $A$ with couplings $2\beta$, external fields $H$ and free boundary condition ($\langle \cdot \rangle_{A^c,2\beta, K}$ is defined analogously); $Z_{A,2\beta,H}$ is its partition function while $Z_{L,n,h}$ is the partition function for $\mu_{L,n,h}$. For the Ising model on a finite graph $G=(V,E)$ with couplings $J_{uv}\geq 0$ for each $uv\in E$ and external fields $g_u\geq 0$ for each $u\in V$, we have
	\[Z_{B \cup C}Z_{B\cap C} \geq Z_B Z_C,\qquad \forall B, C \subset V,\]
	where $Z_B$ is the partition function on the subgraph induced by $B$. This can be proved using the Ahlswede-Daykin four-functions theorem, or one can find a proof in~\cite[Lemma]{Leb74}. 
	
	Therefore, the probability measure $P$ on $2^{V_+}$ satisfies the FKG lattice condition
	\[P(B\cup C) P(B\cap C) \geq P(B) P(C),\qquad \forall B, C \subset V_+.\]
	In particular, this implies that $P$ satisfies the FKG inequality. With the representation~\eqref{eq:stexpectation}, $\langle t_u \rangle_{L,n,h} \geq 0$ because the first GKS inequality \cite{Gri67,KS68} implies that $\langle t_u \rangle_{A^c,2\beta, K}\geq 0$. The second GKS inequality \cite{Gri67,KS68} implies that $\langle t_u \rangle_{A^c,2\beta, K}$ is decreasing in $A$ and $\langle s_v \rangle_{A,2\beta, H}$ is increasing in $A$. So the FKG inequality for $P$ gives
	\begin{align*}
		\langle t_us_v \rangle_{L,n,h}=\sum_{A\subset V_+} P(A) \langle s_v \rangle_{A,2\beta,H} \langle t_u \rangle_{A^c,2\beta,K}&\leq \sum_{A\subset V_+} P(A) \langle s_v \rangle_{A,2\beta,H} \sum_{A\subset V_+} P(A) \langle t_u \rangle_{A^c,2\beta,K}\\
		&=\langle s_v \rangle_{L,n,h} \langle t_u \rangle_{L,n,h}.
	\end{align*}
	
	This implies
	\[\frac{\partial}{\partial H_v}\langle t_u \rangle_{L,n,h}=\langle t_us_v \rangle_{L,n,h}-\langle s_v \rangle_{L,n,h} \langle t_u \rangle_{L,n,h}\leq 0.\]
	The proof of the other inequality in the lemma is similar.
\end{proof}

\begin{lemma}\label{lem:lowerbd}
	Suppose the external field $h\in \R^{\T_L \times I_n}$ satisfies $H_u \geq 0$ and $K_u \geq 0$ for each $u \in V_+$. Let $u,v\in V_+$ with $u\sim v$. Then
	\[\langle t_u \rangle_{L,n,h} \geq \frac{\tanh(2\beta)}{1+\exp[4d\beta+2\|h\|_{\infty}]}\langle t_v \rangle_{L,n,h}.\]
\end{lemma}
\begin{proof}
	Conditioning on all spins except $\sigma_u$ and then averaging the resulting uniform bound give
	\begin{align}\label{eq:uinA^cbd}
		P(u \in A^c)=\mu_{L,n,h}(\sigma_u\neq \sigma_{\bar{u}})&\geq \frac{\exp[-2d\beta-\|h\|_{\infty}]}{\exp[-2d\beta-\|h\|_{\infty}]+\exp[2d\beta+\|h\|_{\infty}]}\nonumber\\
		&=\frac{1}{1+\exp[4d\beta+2\|h\|_{\infty}]}.
	\end{align}
	
	By \eqref{eq:stexpectation} and the second GKS inequality,
	\begin{align*}
		\langle t_u \rangle_{L,n,h}=\sum_{A\subset V_+} P(A) \langle t_u \rangle_{A^c,2\beta,K}
		\geq \sum_{A\subset V_+} P(A) \ind{u\in A^c, v \in A^c}\langle t_u t_v \rangle_{A^c,2\beta,K} \langle t_v \rangle_{A^c,2\beta,K}.
	\end{align*}
	Since $u\sim v$, for every $A$ such that $u, v \in A^c$, the second GKS inequality gives
	\[\langle t_u t_v \rangle_{A^c,2\beta,K}\geq\langle t_u t_v \rangle_{\{u,v\},2\beta,0}=\frac{e^{2\beta}-e^{-2\beta}}{e^{2\beta}+e^{-2\beta}}=\tanh(2\beta). \]
	Plugging this into the last displayed inequality and using the FKG inequality, we get
	\begin{align*}
		\langle t_u \rangle_{L,n,h} &\geq \tanh(2\beta) \sum_{A\subset V_+} P(A) \ind{ u\in A^c} \langle t_v \rangle_{A^c,2\beta,K} \\
		&\geq \tanh(2\beta)  \sum_{A\subset V_+} P(A) \ind{ u\in A^c} \sum_{A\subset V_+} P(A) \langle t_v \rangle_{A^c,2\beta,K}\\
		&=\tanh(2\beta) P(u\in A^c) \langle t_v \rangle_{L,n,h}.
	\end{align*}
	Together with \eqref{eq:uinA^cbd}, this completes the proof of the lemma.
\end{proof}
\begin{proposition}\label{prop:slabandDob}
	Fix $D\geq 3$ and $\beta>0$. For each $l \in \N$, $n>l$ and $u=(x,z) \in \Lambda_l$ with $z\geq 0$,
	\[0 \leq \nu_n(t_u) \leq \langle \sigma_u \rangle^{\mathrm{Dob}}_{\Lambda_l,\beta}.\]
\end{proposition}
\begin{proof}
	Fix $n>l\geq 1$. Choose $L\in \N$ such that $\Lambda_{l+1} \subset \T_L \times I_n$. Recall $h_j$ defined in~\eqref{eq:h_jdef}. In \eqref{eq:mu_Lnhdef}, set $h_{(x,z)}=h_j(z)$ for each $(x,z) \in \T_L\times I_n$.  Then the corresponding $H_u, K_u$ are positive by~\eqref{eq:h_jproperties}. We set
	\[H_u \downarrow 0 \text{ for each }u \in V_+, \ K_u \uparrow \infty \text{ for each }u \in V_+\setminus \Lambda_l.\]
	Lemma~\ref{lem:mon} implies that these operations only increase the expectation of $t_u$; in particular
	\begin{equation}\label{eq:t_uoperations}
		\langle t_u \rangle_{L,n,h'_j} \geq \langle t_u \rangle_{L,n,h_j}, \qquad \forall u \in V_+,
	\end{equation}
	where $h'_j$ denotes the external fields after applying these limits. In the limits, $s_u=0$ and $t_u=1$ for each $u \in V_+\setminus \Lambda_l$. So $\sigma_u=+1=-\sigma_{\bar{u}}$ for each $u \in V_+\setminus \Lambda_l$. Since $n>l$, the antiferromagnetic couplings in~\eqref{eq:mu_Lnhdef} have both endpoints outside $\Lambda_l$. The resulting measure is the Ising measure on $\Lambda_l$ with the Dobrushin boundary condition and external fields $K_u/2$ at each $u=(x,z)\in \Lambda_l$ with $z\geq0$ and $-K_u/2$ at its reflected image $\bar{u}$. This measure (denoted by $\langle \cdot \rangle^{\text{Dob}}_{\Lambda_l,\beta,K}$) is invariant under the reflection with respect to $u_d=-1/2$ followed by spin flip, so
	\begin{equation}\label{eq:t_utoDob}
		\langle t_u \rangle_{L,n,h'_j} = \langle t_u \rangle^{\text{Dob}}_{\Lambda_l,\beta,K}=\langle \sigma_u \rangle^{\text{Dob}}_{\Lambda_l,\beta,K},\qquad \forall u=(x,z) \in \Lambda_l \text{ with }z\geq 0.
	\end{equation}
	
	Now let $L\rightarrow \infty$ along the subsequence defining $\nu_{n,h_j}$, and then let $j\rightarrow \infty$ along the subsequence defining $\nu_n$. Then by \eqref{eq:h_jdef}, $K_u \rightarrow 0$ for each such $u$. Lemma~\ref{lem:mon},~\eqref{eq:t_uoperations} and~\eqref{eq:t_utoDob} complete the proof of the proposition.
\end{proof}

\section{Proof of the main result}\label{sec:proof}
Fix $d\geq 4$ and $\beta>\beta_c(d)$. We first choose a subsequence along which $\nu_n$ converges locally to a measure $\nu$ as $n\rightarrow\infty$. Since each $\nu_n$ is invariant under horizontal translations, so is its local limit $\nu$. 
\begin{claim}
	$\nu$ is a Gibbs state for the zero-field Ising model on $\Z^d$.
\end{claim}
\begin{proof}
	Let $(\pi_{\Lambda})_{\Lambda \Subset \Z^d}$ be the ordinary Ising specification at inverse temperature $\beta$ and no external field. For each finite $\Lambda \subset \Z^d$ and each local function $f$, choose $n$ large enough such that $\Lambda$, its nearest neighbors, and the support of $f$ all lie inside the slab $\Z^D\times I_n$ and away from its antiperiodic boundary. The DLR equation for $\nu_n$ is
	\[\nu_n(f)=\nu_n(\pi_{\Lambda}f)=\int \pi_{\Lambda}(f | \eta) \nu_n(d\eta).\]
	Both $f$ and $\pi_{\Lambda}f$ are local functions, so we get the DLR equation for $\nu$ as $n\rightarrow\infty$ .
\end{proof}
Proposition~\ref{prop:covestimate} implies that
\begin{equation}\label{eq:nuvariancebound}
	\Var_{\nu}\left(\sum_{z} a_z\sigma_{0,z}\right) \leq G_D(0)\sum_{z} a_z^2/\beta,\ \forall a \text{ with finite support},
\end{equation}
\begin{equation}\label{eq:nusummean}
	\sum_{z=-M}^{-1} \nu(\sigma_{0,z}) \leq 0 \leq \sum_{z=0}^{M-1} \nu(\sigma_{0,z}), \ \forall M \geq 1.
\end{equation}

Letting $L\rightarrow \infty$ along the subsequence defining $\nu_{n,h_j}$, then letting $j\rightarrow\infty$ along the subsequence defining $\nu_n$, and finally letting $n\rightarrow\infty$ along the subsequence defining $\nu$, we obtain from Lemmas~\ref{lem:mon} and \ref{lem:lowerbd} that
\begin{equation}\label{eq:nusinglemean}
	\nu(t_{0,z})\geq 0, \  \nu(t_{0,z}) \geq \frac{\tanh(2\beta)}{1+\exp[4d\beta]}\nu(t_{0,z+1}), \ \forall z\geq0.
\end{equation}
\begin{lemma}\label{lem:nupostivemean}
	The state $\nu$ satisfies $\nu(t_{0,0})>0$.
\end{lemma}
\begin{proof}
	Suppose $\nu(t_{0,0})=0$. Then~\eqref{eq:nusinglemean} implies $\nu(t_{0,z})=0$ for each $z\geq 0$. So the two sums in \eqref{eq:nusummean} are equal, and their signs force both sums to vanish. Differences of successive partial sums give
	\begin{equation}\label{eq:numeanzero}
		\nu(\sigma_{0,z})=0,\qquad \forall z \in \Z.
	\end{equation}
	
	Let $\tau_{ke_d} \nu$ be the translation of $\nu$ by $ke_d$. Define
	\[\bar{\nu}_N:=\frac{1}{2N+1}\sum_{k=-N}^N \tau_{ke_d} \nu, \qquad \forall N\in \N.\]
	Let $\bar{\nu}$ be a subsequential local limit of $\bar{\nu}_N$. Since Gibbs states for a finite-range specification form a convex, weakly closed set (see, e.g., ~\cite[ Lemmas 6.27~and~6.28, and Theorem~6.56]{FV18}), $\bar{\nu}$ is a Gibbs state for the zero-field Ising model on $\Z^d$. $\bar{\nu}$ is invariant under horizontal translations because each $\bar{\nu}_N$ is. $\bar{\nu}$ is also invariant under vertical translations since for each local function $f$,
	\[|(\tau_{e_d}\bar{\nu}_N)(f)-\bar{\nu}_N(f)| \leq 2 \|f\|_{\infty} /(2N+1) \rightarrow 0 \text{ as } N\rightarrow \infty.\]
	
	Moreover, $\bar{\nu}(\sigma_{0,0})=0$ by~\eqref{eq:numeanzero}. Therefore, the main result of \cite{Bod06} implies that
	\begin{equation}\label{eq:nudecomposition}
		\bar{\nu}=(\mu^+_{\beta}+\mu^-_{\beta})/2,
	\end{equation}
	where $\mu^+_{\beta}$ and $\mu^-_{\beta}$ are the infinite-volume Ising measures obtained from all plus and all minus boundary conditions respectively.
	
	Define
	\[X_M:=\sum_{z=1}^M \sigma_{0,z}, \qquad \forall M \in \N.\]
	Then $\bar{\nu}(X_M)=0$ since $\bar{\nu}_N(X_M)=0$ for each $N$ by~\eqref{eq:numeanzero}. By~\eqref{eq:numeanzero} and~\eqref{eq:nuvariancebound},
	\[\Var_{\tau_{ke_d}\nu}(X_M)=(\tau_{ke_d} \nu)(X_M^2) \leq MG_D(0)/\beta, \qquad \forall k\in \Z.\]
	
	Thus,
	\[\bar{\nu}_N(X_M^2) \leq M G_D(0)/\beta, \qquad \forall N,M \in \N.\]
	Letting $N\rightarrow\infty$ along the subsequence defining $\bar{\nu}$, we obtain
	\begin{equation}\label{eq:X_Mlowerbound}
		\bar{\nu}(X_M^2) \leq M G_D(0)/\beta, \qquad \forall M\in \N.
	\end{equation}
	
	On the other hand, \eqref{eq:nudecomposition}, the spin-flip invariance of $X_M^2$ and the Cauchy-Schwarz inequality give
	\[\bar{\nu}(X_M^2)=(\mu^+_{\beta}(X_M^2)+\mu^-_{\beta}(X_M^2))/2=\mu^+_{\beta}(X_M^2)\geq[\mu^+_{\beta}(X_M)]^2=M^2[\mu^+_{\beta}(\sigma_{0,0})]^2,\]
	which contradicts \eqref{eq:X_Mlowerbound} for all $M>G_D(0)/(\beta[\mu^+_{\beta}(\sigma_{0,0})]^2)$; here $\mu^+_{\beta}(\sigma_{0,0})>0$ since we assume $\beta>\beta_c(d)$.  This contradiction implies that $\nu(t_{0,0})>0$.
\end{proof}

We are ready to prove our main result.
\begin{proof}[Proof of Theorem~\ref{thm:main}]
	Proposition~\ref{prop:slabandDob} gives
	\[\langle \sigma_{0,0} \rangle^{\text{Dob}}_{\Lambda_l,\beta} \geq \nu_n(t_{0,0}), \qquad \forall n>l\geq 1.\]
	Letting $n\rightarrow\infty$ along the subsequence defining $\nu$ and then using Lemma~\ref{lem:nupostivemean}, we get
	\[\langle \sigma_{0,0} \rangle^{\text{Dob}}_{\Lambda_l,\beta} \geq \nu(t_{0,0})>0, \qquad \forall l\in \N.\]
	Since $\mu^{\text{Dob}}_{\Lambda_l,\beta}$ is invariant under the reflection with respect to $u_d=-1/2$ followed by a spin-flip,
	\begin{equation}\label{eq:Dobat0and-1}
		\langle \sigma_{0,0} \rangle^{\text{Dob}}_{\Lambda_l,\beta}=-\langle \sigma_{0,-1} \rangle^{\text{Dob}}_{\Lambda_l,\beta}\geq \nu(t_{0,0})>0, \qquad \forall l\in \N.
	\end{equation}
	Let $\mu^{\text{Dob}}_{\beta}$ be a subsequential local limit of $\mu^{\text{Dob}}_{\Lambda_l,\beta}$ as $l\rightarrow\infty$. Then $\mu^{\text{Dob}}_{\beta}$  is a Gibbs state; see, e.g., \cite[Section 6.2.1]{FV18}. By~\eqref{eq:Dobat0and-1},
	\[\mu^{\text{Dob}}_{\beta}(\sigma_{0,0})=-\mu^{\text{Dob}}_{\beta}(\sigma_{0,-1})>0.\]
	Therefore, $\mu^{\text{Dob}}_{\beta}$  is not invariant under vertical translations.
\end{proof}

\subsection*{AI statement}
The overall idea of using correlation inequalities and the infrared bound to prove the main result was due to the authors. The crucial idea of using an antiperiodic slab was obtained with the help of AI tools. The entire manuscript was written by the authors. AI tools were used for proofreading.

\section*{Acknowledgments}
This research was partially supported by the National Natural Science Foundation of China (No. 12271284 and No. 12226001). 

\bibliographystyle{abbrv}
\bibliography{references}

\end{document}